\documentclass[11pt,reqno]{amsart}

\usepackage{amsmath, amsfonts, amsthm, amssymb, graphicx}

\usepackage{amssymb}
\usepackage{amsmath}
\usepackage{color}
\usepackage{graphics}
\usepackage{epsfig}
\usepackage{hyperref}
\usepackage{epstopdf}

\newtheorem{claim}{\bf \t}[part]

\theoremstyle{plain}
\newtheorem{Theorem}{Theorem}[section]
\newtheorem{Lemma}{Lemma}[section]
\newtheorem{Proposition}{Proposition}[section]

\theoremstyle{definition}
\newtheorem{Definition}{Definition}[section]

\theoremstyle{remark}
\newtheorem{Remark}{Remark}[section]

\theoremstyle{definition}
\newtheorem{Example}{Example}[section]

\numberwithin{equation}{section}

\def\R{\mathbb{R}}

\begin{document}

\title[Isometric Immersions of Surfaces]
{Isometric Immersions of Surfaces with Two Classes of Metrics and Negative Gauss Curvature}
\author{Wentao Cao \and Feimin Huang \and Dehua Wang}
\address{W. Cao, Institute of Applied Mathematics, AMSS, CAS, Beijing 100190, China.}
\email{cwt@amss.ac.cn}

\address{F. Huang, Institute of Applied Mathematics, AMSS, CAS, Beijing 100190, China.}
 \email{\tt fhuang@amt.ac.cn}

\address{D. Wang, Department of Mathematics, University of Pittsburgh,
                Pittsburgh, PA 15260, USA.}
\email{\tt dwang@math.pitt.edu}


\date{}
\thanks{}

\begin{abstract}
The isometric immersion of two-dimensional Riemannian manifolds or surfaces with negative Gauss curvature into the three-dimensional Euclidean space is studied in this paper. The global weak solutions to the Gauss-Codazzi equations with large data in $L^\infty$ are obtained through the vanishing viscosity method and the compensated compactness framework. The $L^\infty$ uniform estimate and $H^{-1}$ compactness are established through a transformation of state variables and construction of proper invariant regions for two types of given metrics including the catenoid type and the helicoid type. The global weak solutions in $L^\infty$ to the Gauss-Codazzi equations yield the $C^{1,1}$ isometric immersions of  surfaces with the given metrics.
\end{abstract}

\maketitle

\section{Introduction } \label{S1}

Isometric embedding or immersion into
$\mathbb{R}^3$ of  two-dimensional Riemannian manifolds (or surfaces) $\mathcal{M}^2$
is a well known classical problem in differential geometry (cf.  Cartan \cite{Cartan} in 1927, Codazzi \cite{Codazzi} in 1860, Janet \cite{Janet} in 1926, Mainardi \cite{Mainardi} in 1856, Peterson \cite{Peterson} in 1853), with emerging applications in shell theory, computer graphics, biological leaf growth and  protein folding in biology and so on (cf. \cite{Ghrist1,VM}).
The classical surface theory indicates that for the given metric,  the isometric embedding or immersion can be realized if the first fundamental form and the second fundamental form satisfy the Gauss-Codazzi equations (cf. \cite{BS,M1,M2,PS}).
There have been many  results  for the isometric embedding of
surfaces with positive Gauss curvature, which can be studied
by solving an elliptic  problem of the Darboux equation or the
Gauss-Codazzi equations;  see \cite{HH} and the references therein.
When the Gauss curvature is negative or changes signs, there are only a few studies in literature.
The case where the Gauss curvature is negative can be formulated into a hyperbolic
 problem of nonlinear partial differential equations, and the case in which the Gauss
curvature changes signs becomes solving nonlinear partial differential
equations of mixed elliptic-hyperbolic type. Han in \cite{Han} obtained local isometric
embedding of surfaces with Gauss curvature changing sign cleanly.
Hong in \cite{H}  proved  the  isometric immersion in $\mathbb{R}^3$ of completely
negative curved surfaces with the negative Gauss curvature decaying at a certain rate
in the time-like direction, and the  $C^1$ norm of initial data is  small so
that he can obtain the  smooth solution.  Recently, Chen-Slemrod-Wang in \cite{CSW} developed
a general method, which combines a fluid dynamic formulation of conservation
laws for the Gauss-Codazzi system with a compensated compactness framework,
to realize the isometric immersions in $\mathbb{R}^3$ with negative Gauss curvature.
Christoforou in \cite{Christoforou} obtained the small BV solution to the Gauss-Codazzi system with the same catenoid type metric as in \cite{CSW}.
See \cite{Dong, Efimov-a, Efimov, H, CSLin, PS, Roz, Yau00} for other related results on surface embeddings.
For the higher dimensional isometric embeddings we refer the readers to \cite{BBG1, BBG2, BGY, CSW2, CSW3, Eisenhart, Gromov86, NM, NM2, Nash56, Poole} and the references therein.

Recall that in  Chen-Slemrod-Wang \cite{CSW},  the  $C^{1,1}$ isometric immersion  of surfaces was obtained  for the  catenoid type metric:
$$ds^2=E(y)dx^2+cE(y)dy^2$$
 and negative Gauss curvature $$K(y)=-k_0E(y)^{-{\beta}^2}$$
  with index $\beta>\sqrt{2}$, $c=1$ and  $k_0>0$.
The goal of this paper is to study the isometric immersion of surfaces with negative Gauss curvature and more general metrics.
To this end, we shall apply directly the artificial vanishing viscosity method (\cite{Dafermos-book}) and introduce a  transformation of state variables (see e.g. \cite{HH}).  The vanishing viscosity limit will be obtained through the compensated compactness method (\cite{Ball, Chen1, DCL, DiPerna1, DiPerna3, evans, huangwang, LPS, LPT, Mu2, Ta1}). We shall establish the $L^\infty$ estimate of the viscous approximate solutions by studying the Riemann invariants to find the invariant region (see e.g \cite{Smoller-book})  in the new state variables. Then we prove the $H^{-1}$ compactness of the viscous approximate solutions, and finally apply the compensated compactness framework in \cite{CSW} to obtain the weak solution of the Gauss-Codazzi equations,  and thus realize the isometric immersion of surfaces into $\mathbb{R}^3$.
In the new state variables, using the vanishing artificial viscosity method, we are able to obtain the weak solution of the Gauss-Codazzi system for two  classes of metrics,   
that is, for the catenoid type metric with  $\beta\geq\sqrt{2}$ and $c>0$, and for the  helicoid type metric
$$ds^2=E(y)dx^2+dy^2$$
and $$K(y)=-k_0E(y)^a$$ with $a\leq-2$ and $k_0>0$.
An important step to achieve this development is to find the invariant regions of the viscous approximate solutions  for the wider classes of metrics, for which the introduction of new state variables $(u,v)$ plays a crucial role.
We note that the $L^\infty$ solution of the Gauss-Codazzi equations
for the given  metric in $C^{1,1}$ yields the  $C^{1,1}$
isometric immersion from  the fundamental theorem of surfaces
 by Mardare \cite{M1,M2}.

The paper is organized as follows.
In Section 2, we introduce a new formulation of the
Gauss-Codazzi system and provide  the viscous approximate solutions. In Section 3,
we establish the  $L^\infty$ estimate for the two types of metrics
to get the global existence of the equations with viscous terms. In Section 4,
we prove the  $H^{-1}_{loc}$ compactness. In Section 5, combining
the above two estimates and the compensated compactness framework we state and prove our theorems on the existence of weak solutions for the surfaces  with the catenoid type and helicoid type metrics.

\bigskip

\section{Reformulation of the Gauss-Codazzi System}\label{S2}

As in \cite{HH} and \cite{CSW}, the isometric embedding problem
for the two-dimensional Riemannian manifolds (or surfaces) in $\mathbb{R}^3$ can be
formulated through the Gauss-Codazzi system
using the fundamental theorem of surface theory (cf.  Mardare \cite{M1,M2}).

Let $\Omega\subset\mathbb{R}^2$ be an open set and  $(x,y)\in\Omega$.
For a two-dimensional surface defined on $\Omega$ with the given metric in the first fundamental form:
\begin{equation*}
I=Edx^2+2Fdxdy+Gdy^2,
\end{equation*}
where $E,F,G$ are differentiable functions of $(x,y)$ in $\Omega$.
If  the second fundamental form of the surface is
\begin{equation*}
I\!\!I=h_{11}dx^2+2h_{12}dxdy+h_{22}dy^2,
\end{equation*}
where $h_{11}, h_{12}=h_{21}, h_{22}$ are also functions of $(x,y)$ in $\Omega$,
then the Gauss-Codazzi system has the following form:
\begin{equation}\label{Gc}
\begin{cases}
 &\displaystyle M_x-L_y=\Gamma^2_{22}L-2\Gamma^2_{12}M+\Gamma^2_{11}N,\\
 &\displaystyle N_x-M_y=-\Gamma^1_{22}L+2\Gamma^1_{12}M-\Gamma^1_{11}N,\\
 &\displaystyle LN-M^2=K.
\end{cases}
\end{equation}
Here and in the rest of the paper $\Box_x$,  $\Box_y$ stand for the partial
 derivatives of corresponding function $\Box$ with respect to $x, y$ respectively.
In \eqref{Gc},
\begin{equation*}
L=\frac{h_{11}}{\sqrt{|g|}},\quad M=\frac{h_{12}}{\sqrt{|g|}}, \quad N=\frac{h_{22}}{\sqrt{|g|}},
\end{equation*}
$$|g|=EG-F^2>0,$$
 $K=K(x,y)$ is the Gauss curvature which is determined
by $E, F, G$ according to the Gauss's Theorem  Egregium (\cite{doC1992,Roz}),
 $\Gamma^k_{ij} \; (i, j, k=1,2)$ are the
Christoffel symbols 
given by the following formulas (\cite{HH}):
\begin{equation}\label{Chr1}
\begin{split}
& \Gamma^1_{11}=\frac{GE_x-2FF_x+FE_y}{2(EG-F^2)},
 \quad \Gamma^2_{11}=\frac{2EF_x-EE_y-FE_x}{2(EG-F^2)},\\
& \Gamma^1_{12}=\Gamma^1_{21}=\frac{GE_y-FG_x}{2(EG-F^2)},
 \quad \Gamma^2_{12}=\Gamma^2_{21}=\frac{EG_x-FE_y}{2(EG-F^2)}, \\
& \Gamma^1_{22}=\frac{2GF_y-GG_x-FG_x}{2(EG-F^2)},
 \quad \Gamma^2_{22}=\frac{EG_y-2FF_y+FG_x}{2(EG-F^2)}.
 \end{split}
\end{equation}
As in Mardare \cite{M1,M2} and Chen-Slemrod-Wang \cite{CSW}, the fundamental theorem of surface theory holds when $(h_{ij})$ or $L, M, N$ are  in $L^{\infty}$
for the given $E,F,G\in C^{1,1}$,  and  the immersion surface
is $C^{1,1}$. Therefore it suffices to find the solutions $L, M, N$ in
$L^{\infty}$ to the Gauss-Codazzi system to realize the
surface with the given metric.

In this paper, we  consider the isometric immersion into $\mathbb{R}^3$ of a
two-dimensional Riemanian manifold  with negative Gauss curvature.
We write  the negative Gauss curvature  as
\begin{equation}\label{K}
K=-\gamma^2
\end{equation}
 with $\gamma>0$ and $\gamma\in C^{1,1}(\Omega)$,
and rescale $L, M, N$ as
\begin{equation}
\tilde{L}=\frac{L}{\gamma},\quad \tilde{M}=\frac{M}{\gamma},\quad \tilde{N}=\frac{N}{\gamma}.
\end{equation}
Then the third equation (or the Gauss equation) of system (\ref{Gc}) becomes
\begin{equation}\label{Sta}
\tilde{L}\tilde{N}-\tilde{M}^2=-1.
\end{equation}
The other two equations (or the Codazzi equations) of system \eqref{Gc} become
\begin{equation}\label{gc}
\begin{split}
& \tilde{M}_x-\tilde{L}_y=\tilde{\Gamma}^2_{22}\tilde{L}-2\tilde{\Gamma}^2_{12}\tilde{M}+\tilde{\Gamma}^2_{11}\tilde{N},\\
& \tilde{N}_x-\tilde{M}_y=-\tilde{\Gamma}^1_{22}\tilde{L}+2\tilde{\Gamma}^1_{12}\tilde{M}-\tilde{\Gamma}^1_{11}\tilde{N},
\end{split}
\end{equation}
where
\begin{equation}\label{Chr2}
\begin{split}
\tilde{\Gamma}^2_{22}=\Gamma^2_{22}+\frac{\gamma_y}{\gamma},\quad \tilde{\Gamma}^2_{12}=\Gamma^2_{12}+\frac{\gamma_x}{2\gamma},\quad \tilde{\Gamma}^2_{11}=\Gamma^2_{11},\\
\tilde{\Gamma}^1_{22}=\Gamma^1_{22},\quad \tilde{\Gamma}^1_{12}=\Gamma^1_{12}+\frac{\gamma_y}{2\gamma},\quad \tilde{\Gamma}^1_{11}=\Gamma^1_{11}+\frac{\gamma_x}{\gamma}.
\end{split}
\end{equation}

Consider the following viscous approximation  of system \eqref{Sta}-\eqref{gc}  with artificial viscosity:
\begin{equation}\label{gcv}
\begin{cases}
& \tilde{L}_y-\tilde{M}_x=\varepsilon\tilde{L}_{xx}
  -\tilde{\Gamma}^2_{22}\tilde{L}+2\tilde{\Gamma}^2_{12}\tilde{M}-\tilde{\Gamma}^2_{11}\tilde{N},\\
& \tilde{M}_y-\tilde{N}_x=\varepsilon\tilde{M}_{xx}
 +\tilde{\Gamma}^1_{22}\tilde{L}-2\tilde{\Gamma}^1_{12}\tilde{M}+\tilde{\Gamma}^1_{11}\tilde{N},\\
 &\tilde{L}\tilde{N}-\tilde{M}^2=-1,
\end{cases}
\end{equation}
where $\varepsilon>0$.
Our goal is to apply the vanishing viscosity method to the smooth solutions of \eqref{gcv} to obtain the $L^\infty$ solution of
\eqref{Sta}-\eqref{gc}.
The eigenvalues  of the system  \eqref{gc}  for $\tilde{L}\neq 0$ are
$$\lambda_{\pm}=\frac{-\tilde{M}\pm1}{\tilde{L}}, $$
and the right eigenvectors are
$$\vec{r}_{\pm}=(1,-\lambda_{\pm})^{T}.$$
A direct calculation shows
$$\bigtriangledown \lambda_{\pm}\cdot\vec{r}_{\pm}=
\frac{-\tilde{M}\pm1}{-\tilde{L}^2}\cdot1+
\frac{-1}{\tilde{L}}\cdot\frac{-\tilde{M}\pm1}{-\tilde{L}}=0,$$
thus we can take $\lambda_\pm$ as the Riemann invariants.

Introduce the new variables:
\begin{equation}
u:=-\frac{\tilde{M}}{\tilde{L}},\quad v:=\frac{1}{\tilde{L}},
\end{equation}
 then
\begin{equation}
\tilde{L}=\frac{1}{v},\quad \tilde{M}=\frac{-u}{v},\quad \tilde{N}=\frac{u^2-v^2}{v}.
\end{equation}
Thus
\begin{equation}\label{1.1}
\begin{split}
&\tilde{L}_y=-\frac{v_y}{v^2}, \qquad \tilde{L}_{xx}=\frac{2v_x^2-vv_{xx}}{v^3},\\
&\tilde{N}_x=\frac{(2uu_x-2vv_x)v-(u^2-v^2)v_x}{v^2}=\frac{2uvu_x-v^2v_x-u^2v_x}{v^2},
\end{split}
\end{equation}
and
\begin{equation}\label{1.2}
\begin{split}
&\tilde{M}_y=\frac{uv_y-vu_y}{v^2}, \quad \tilde{M}_x=\frac{uv_x-vu_x}{v^2},\\
&\tilde{M}_{xx}=\frac{uv_{xx}-vu_{xx}}{v^2}-\frac{2uv_x^2-2vu_xv_x}{v^3}.
\end{split}
\end{equation}
Substituting (\ref{1.1}) and (\ref{1.2}) into the system (\ref{gcv}), we get
\begin{equation}\label{2.1}
\frac{-v_y}{v^2}-\frac{uv_x-vu_x}{v^2}=\frac{2\varepsilon v_x^2-\varepsilon vv_{xx}}{v^3}
-\tilde{\Gamma}^2_{22}\frac{1}{v}-2\tilde{\Gamma}^2_{12}\frac{u}{v}-\tilde{\Gamma}^2_{11}\frac{u^2-v^2}{v},
\end{equation}
and
\begin{equation}\label{2.2}
\begin{split}
&\frac{uv_y-vu_y}{v^2}-\frac{2uvu_x-v^2v_x-u^2v_x}{v^2}\\
&=\frac{\varepsilon uv_{xx}-\varepsilon vu_{xx}}{v^2}-
\frac{2\varepsilon uv_x^2-2\varepsilon vu_xv_x}{v^3}
+\tilde{\Gamma}^1_{22}\frac{1}{v}+2\tilde{\Gamma}^1_{12}\frac{u}{v}+\tilde{\Gamma}^1_{11}\frac{u^2-v^2}{v}.
\end{split}
\end{equation}
Multiplying (\ref{2.1}) by $-v^2$,  we get
\begin{equation}\label{3.1}
v_y+uv_x-vu_x=\varepsilon v_{xx}-\frac{2\varepsilon v_x^2}{v}+
\tilde{\Gamma}^2_{22}v+2\tilde{\Gamma}^2_{12}uv+\tilde{\Gamma}^2_{11}(u^2-v^2)v,
\end{equation}
and multiplying  (\ref{2.2}) by $-v$, one has
\begin{equation}\label{3.2}
\begin{split}
&-\frac{uv_y}{v}+u_y+2uu_x-vv_x-\frac{u^2v_x}{v}\\
&=\varepsilon u_{xx}-\frac{\varepsilon uv_{xx}}{v}
-\frac{2\varepsilon u_xv_x}{v}+\frac{2\varepsilon uv_x^2}{v^2}
 -\tilde{\Gamma}^1_{22}-2\tilde{\Gamma}^1_{12}u-\tilde{\Gamma}^1_{11}(u^2-v^2).
 \end{split}
\end{equation}
Substituting  (\ref{3.1}) into (\ref{3.2}), we obtain
\begin{equation*}
\begin{split}
u_y+uu_x-vv_x=\, &\varepsilon u_{xx}-\frac{2\varepsilon v_xu_x}{v}
-\tilde{\Gamma}^1_{22}-2\tilde{\Gamma}^1_{12}u-\tilde{\Gamma}^1_{11}(u^2-v^2)\\
 &+\tilde{\Gamma}^2_{22}u+2\tilde{\Gamma}^2_{12}u^2+\tilde{\Gamma}^2_{11}(u^2-v^2)u\\
=\, &\varepsilon u_{xx}-\frac{2\varepsilon v_xu_x}{v}
-\tilde{\Gamma}^1_{22}+(\tilde{\Gamma}^2_{22}-2\tilde{\Gamma}^1_{12})u+
(2\tilde{\Gamma}^2_{12}-\tilde{\Gamma}^1_{11})u^2\\
 &+\tilde{\Gamma}^1_{11}v^2+\tilde{\Gamma}^2_{11}(u^2-v^2)u.
\end{split}
\end{equation*}
Therefore we have the following system in the variables $(u,v)$:
\begin{equation}\label{4}
\begin{split}
u_y+(uu_x-vv_x)=&\varepsilon u_{xx}-\frac{2\varepsilon v_xu_x}{v}
  -\tilde{\Gamma}^1_{22}+(\tilde{\Gamma}^2_{22}-2\tilde{\Gamma}^1_{12})u\\
   &+(2\tilde{\Gamma}^2_{12}-\tilde{\Gamma}^1_{11})u^2 +
\tilde{\Gamma}^1_{11}v^2+\tilde{\Gamma}^2_{11}(u^2-v^2)u,\\
v_y+(uv_x-vu_x)=&\varepsilon v_{xx}-\frac{2\varepsilon v_x^2}{v}+
\tilde{\Gamma}^2_{22}v+2\tilde{\Gamma}^2_{12}uv+\tilde{\Gamma}^2_{11}(u^2-v^2)v.
\end{split}
\end{equation}
We note that the $L^\infty$ estimate for  the solutions $\tilde{L}, \tilde{M}$ and $\tilde{N}$  to \eqref{gcv} with $\tilde{L}>0$ is equivalent to  the $L^\infty$ estimate for the solutions $u$ and $v$ to \eqref{4} with $v>0$.
Set
\begin{equation}\label{5}
\begin{split}
&f(u,v)=-\tilde{\Gamma}^1_{22}+(\tilde{\Gamma}^2_{22}-2\tilde{\Gamma}^1_{12})u+
(2\tilde{\Gamma}^2_{12}-\tilde{\Gamma}^1_{11})u^2+\tilde{\Gamma}^1_{11}v^2
+\tilde{\Gamma}^2_{11}(u^2-v^2)u,\\
&g(u,v)=\tilde{\Gamma}^2_{22}v+2\tilde{\Gamma}^2_{12}uv+\tilde{\Gamma}^2_{11}(u^2-v^2)v,
\end{split}
\end{equation}
then the system \eqref{4} becomes
\begin{equation}\label{6}
\begin{cases}
&\displaystyle u_y+(uu_x-vv_x)=f(u,v)+\varepsilon u_{xx}-\frac{2\varepsilon u_xv_x}{v},\\[10pt]
&\displaystyle v_y+(uv_x-vu_x)=g(u,v)+\varepsilon v_{xx}-\frac{2\varepsilon v_x^2}{v}.
\end{cases}
\end{equation}

The local existence of (\ref{gcv}) is standard, and the global existence can be proved if  the $L^\infty$ boundedness of
$\tilde{L},\tilde{M}$ and $\tilde{N}$ (or equivalently the $L^\infty$ boundedness of $u$ and $v$) is established.
The $L^\infty$ uniform bound  will be established  in the next Section \ref{S3}.
The global existence of solution $(\tilde{L},\tilde{M})$ to \eqref{gcv} is equivalent to the global existence of solution $(u, v)$  to \eqref{6} which will be proved in Section \ref{S5}.

\bigskip

\section{$L^{\infty}$  Uniform Estimate} \label{S3}

In order to establish  the $L^\infty$ bound of $u$ and $v$, we need to derive the equations of the Riemann invariants of (\ref{6}). First we  rewrite the system (\ref{6}) in the following form:
\begin{equation}\label{7}
\begin{bmatrix} u_y\\ v_y \end{bmatrix}
+\begin{bmatrix} u & -v\\ -v & u \end{bmatrix}\begin{bmatrix} u_x \\ v_x\end{bmatrix}
=\begin{bmatrix} f(u,v)\\ g(u,v) \end{bmatrix}+\begin{bmatrix} \varepsilon u_{xx}-\frac{2\varepsilon u_xv_x}{v}\\
\varepsilon v_{xx}-\frac{2\varepsilon v_x^2}{v} \end{bmatrix}.
\end{equation}
%
The eigenvalues of \eqref{7} are
\begin{equation*}
\lambda_1=u-v,\quad \lambda_2=u+v,
\end{equation*}
and the Riemann invariants are
$$w=u+v, \quad z=u-v.$$
 Multiply (\ref{7}) by $(w_u, w_v)$ or $(z_u, z_v)$ to obtain the equations satisfied by the Riemann invariants:
\begin{equation}\label{0}
\begin{split}
&w_y+\lambda_1w_x=\varepsilon w_{xx}-\frac{2\varepsilon v_xw_{x}}{v}+f(u,v)+g(u,v),\\
&z_y+\lambda_2z_x=\varepsilon z_{xx}-\frac{2\varepsilon v_xz_{x}}{v}+f(u,v)-g(u,v).
\end{split}
\end{equation}
Then at the critical points of $w$, the first equation of (\ref{0}) becomes
\begin{equation*}
\varepsilon w_{xx}+f(u,v)+g(u,v)=0,
\end{equation*}
and at the critical points of $z$, the second equation of (\ref{0}) becomes
\begin{equation*}
\varepsilon z_{xx}+f(u,v)-g(u,v)=0.
\end{equation*}
Hence, by the parabolic maximum principle (see \cite{LC,PW}), we have
\begin{itemize}
\item[(a)]$w$ has no internal maximum when $f(u,v)+g(u,v)<0,$
\item[(b)]$w$ has no internal minimum when $f(u,v)+g(u,v)>0,$
\item[(c)]$z$ has no internal maximum when $f(u,v)-g(u,v)<0,$
\item[(d)]$z$ has no internal minimum when $f(u,v)-g(u,v)>0.$
\end{itemize}
In order to find the invariant region of $w, z$, we need to analyze the source terms in \eqref{0}, that is,  the signs of $f(u,v)+g(u,v)$ and $f(u,v)-g(u,v)$. We shall consider two types of surfaces that have special metrics with $F\equiv0.$

\subsection{Catenoid type surfaces: $G(y)=cE(y), F=0$, $c>0$} \label{catenoid}
For the surfaces with metrics of the form:
\begin{equation}\label{catenoidmetric}
G(y)=cE(y), \quad F=0
\end{equation}
 with constant $c>0$,
from the formulas of $\Gamma^k_{ij}$ in (\ref{Chr1}) and \eqref{Chr2}, we can easily calculate that
\begin{equation*}
\begin{split}
&\Gamma^2_{22}=\frac{E^\prime}{2E},\qquad  \Gamma^2_{12}=0, \qquad \Gamma^2_{11}=\frac{-E^\prime}{2cE},\\
&\Gamma^1_{22}=0,\qquad  \Gamma^1_{12}=\frac{E^\prime}{2E},\qquad  \Gamma^1_{11}=0,
\end{split}
\end{equation*}
and thus
\begin{equation*}\label{8}
\begin{split}
&\tilde{\Gamma}^2_{22}=\frac{E^\prime}{2E}+\frac{\gamma^\prime}{\gamma}, \qquad \tilde{\Gamma}^2_{12}=0, \qquad \tilde{\Gamma}^2_{11}=\frac{-E^\prime}{2cE},\\
&\tilde{\Gamma}^1_{22}=0, \qquad \tilde{\Gamma}^1_{12}=\frac{E^\prime}{2E}+\frac{\gamma^\prime}{2\gamma}, \qquad \tilde{\Gamma}^1_{11}=0.
\end{split}
\end{equation*}
If we assume  $$\frac{\gamma^\prime}{\gamma}=\alpha \frac{E^\prime}{E},$$
 where $\alpha$ is constant, then, by \eqref{K},  $$2\alpha \frac{E^\prime}{E}=\frac{K^\prime}{K},$$
  that is
  \begin{equation}\label{cm1}
  K(y)=-k_0E(y)^{2\alpha},
  \end{equation}
   with some constant $k_0>0$.
From  the expressions of $f(u,v)$, $g(u,v)$ in \eqref{5}, one has
\begin{equation*}
\begin{split}
&f(u,v)=\frac{-E^\prime}{2E}u-\frac{E^\prime}{2cE}(u^2-v^2)u,\\
&g(u,v)=\frac{E^\prime}{E}(\alpha+\frac{1}{2})v-\frac{E^\prime}{2cE}(u^2-v^2)v.
\end{split}
\end{equation*}
and then
\begin{equation}
\begin{split}
f(u,v)+g(u,v)=-\frac{E^\prime}{2cE}\left(cu-c(2\alpha+1)v+(u^2-v^2)(u+v)\right),\\
f(u,v)-g(u,v)=-\frac{E^\prime}{2cE}\left(cu+c(2\alpha+1)v+(u^2-v^2)(u-v)\right).
\end{split}
\end{equation}
Set
\begin{equation}
\begin{split}
\varphi_1(u,v)=cu-c(2\alpha+1)v+(u^2-v^2)(u+v), \\ \varphi_2(u,v)=cu+c(2\alpha+1)v+(u^2-v^2)(u-v).
\end{split}
\end{equation}
In particular, when $\alpha=-1$,
\begin{equation}
\begin{split}
\varphi_1(u,v)=(u^2-v^2+c)(u+v), \\ \varphi_2(u,v)=(u^2-v^2+c)(u-v).
\end{split}
\end{equation}
If we assume
\begin{equation}\label{cm2}
\frac{E'}{E}<0,
\end{equation}
 then the signs of $$f(u,v)+g(u,v), \quad f(u,v)-g(u,v)$$
 depend only on $\varphi_1(u,v),\varphi_2(u,v)$ respectively.

We now derive and sketch the invariant regions.

When $\alpha=-1$, first we can easily sketch the level sets of $w,z$ in the $u-v$ plane in
Figure 1.
\begin{figure}
\begin{center}
  \includegraphics[width=16cm,height=8cm]{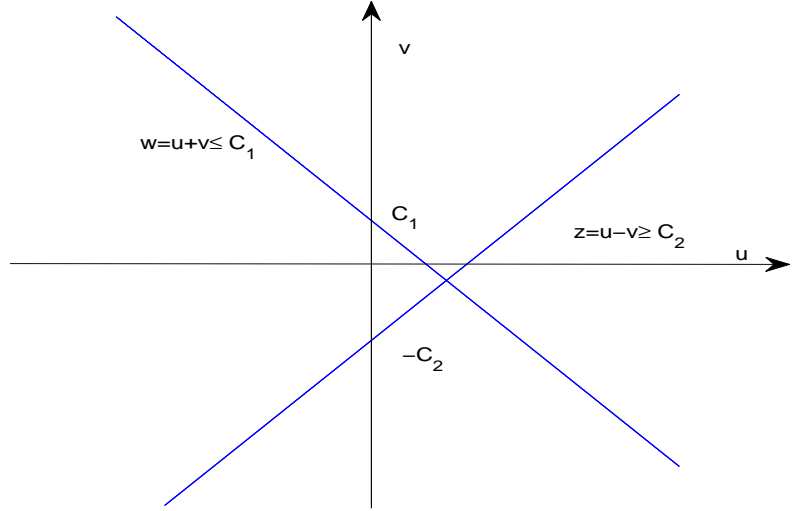}
  \caption{Graphs of level sets of $w$ and $z$}\label{figure1}
  \end{center}
\end{figure}
\begin{figure}
\begin{center}
  \includegraphics[width=16cm,height=8cm]{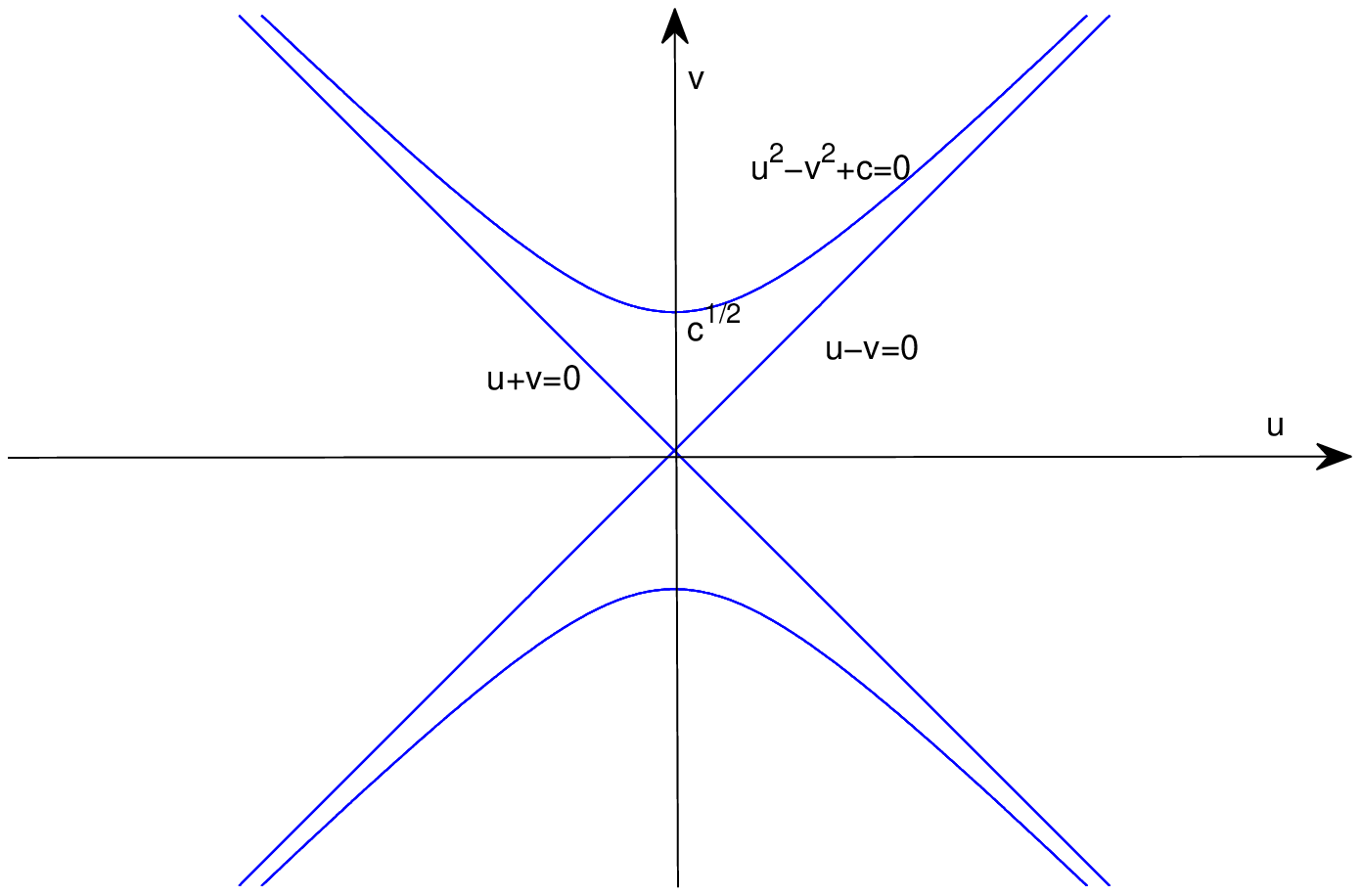}
  \caption{Graphs of $\varphi_1(u,v)=0, \varphi_2(u,v)=0$ when $\alpha=-1$}\label{figure2}
  \end{center}
\end{figure}
\begin{figure}
\begin{center}
  \includegraphics[width=16cm,height=8cm]{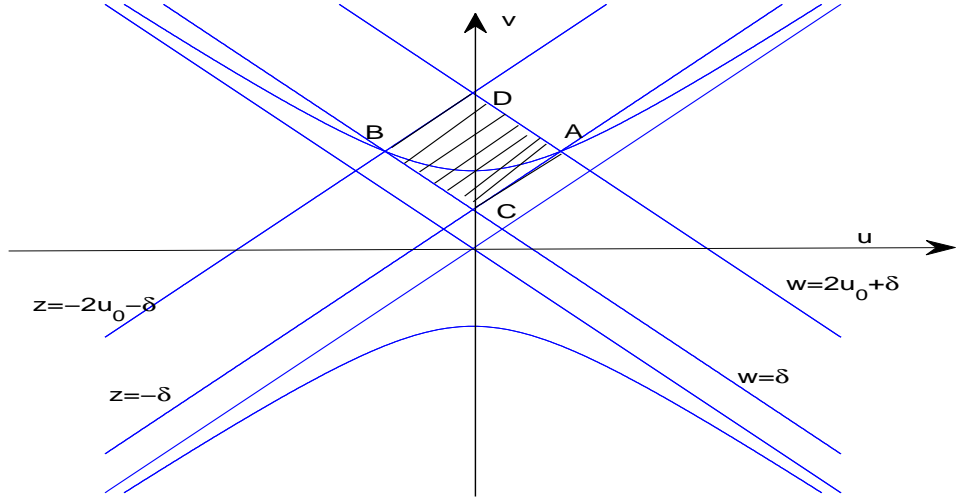}\\
  \caption{Invariant region when $\alpha=-1$}\label{figure3}
  \end{center}
\end{figure}
To obtain the invariant region we need to sketch the graphs of
$$\varphi_1(u,v)=0, \, \varphi_2(u,v)=0$$ in the $u-v$ plane in Figure 2.
Now let us find the invariant region in the upper half-plane.  We
draw a straight line parallel to $u-v=0$ passing through the point $$C=(0,\delta)$$ with
$0<\delta<\sqrt{c}$ intersecting the hyperbola at the point
$$A=(u_0, u_0+\delta),$$
and  similarly, we get the point
$$B=(-u_0, u_0+\delta),$$
where
\begin{equation}\label{u01}
u_0=u_0(-1)=\frac{c-\delta^2}{2\delta}.
\end{equation}
  We then draw a straight
line perpendicular to $u-v=0$  passing through the point $A$ and another straight line
perpendicular to $u+v=0$ through $B$. Then the two lines
intersect the $v-$axis at the same point $$D=\left(0, \frac{c-\delta^2}{2\delta}\right);$$
 see Figure 3.
We  see that the square $ACBD$ in Figure 3 is an  invariant
region. Therefore we get the $L^\infty$ estimate of $(u,v)$, that is,
\begin{equation}
-\frac{c-\delta^2}{2\delta}\leq u \leq \frac{c-\delta^2}{2\delta},
\quad \delta \leq v \leq \frac{c}{\delta}.
\end{equation}

When $\alpha<-1$,
the graphs of $\varphi_1(u,v)=0, \, \varphi_2(u,v)=0$ in the $u-v$
plane look like the curves marked with $1$ and $2$ respectively in Figure 4. Similarly to the case $\alpha=-1$, we
draw a straight line parallel to $u-v=0$ passing through the point $$C=(0,\delta)$$ with
$0<\delta<\sqrt{-2c\alpha-c}$ intersecting the curve
$\varphi_1(u,v)=0$ at the point $$A=(u_0, u_0+\delta).$$  Then we draw a straight
line parallel to $u+v=0$ passing through the point $C=(0,\delta)$ intersecting the curve
$\varphi_2(u,v)=0$ at the point $$B=(-u_0, u_0+\delta).$$  Moreover, we draw a
straight line perpendicular to $u-v=0$ passing through point $A$ and another
straight line perpendicular to $u+v=0$ passing through point $B$. Then the two
lines   intersect the $v-$axis at the same point $$D=(0, 2u_0+\delta).$$
Thus, the square $ACBD$ in Figure \ref{figure5} is an invariant region,
and  now the $L^\infty$ estimate of $(u,v)$ is
\begin{equation}\label{u0}
-u_0\leq u\leq u_0, \quad \delta \leq v\leq 2u_0+\delta,
\end{equation}
where $$0<\delta<\sqrt{-2c\alpha-c},$$
\begin{equation}\label{u02}
u_0=u_0(\alpha)=-\frac{c\alpha+2\delta^2-\sqrt{c^2\alpha^2-4c\alpha\delta^2-4c\delta^2}}{4\delta}.
\end{equation}
\begin{figure}
\begin{center}
  \includegraphics[width=15cm,height=8cm]{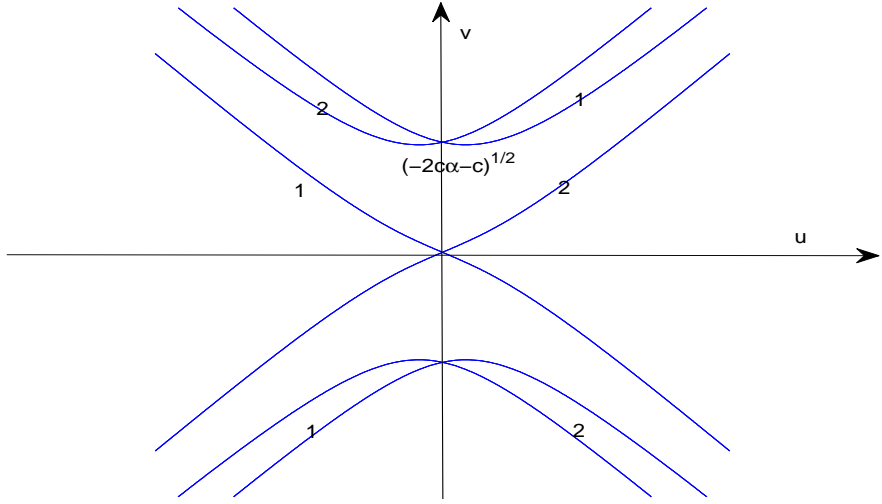}\\
  \caption{Graphs of $\varphi_1(u,v)=0, \varphi_2(u,v)=0$ when $\alpha<-1$.}
  \label{figure4}
  \end{center}
\end{figure}
\begin{figure}
\begin{center}
  \includegraphics[width=15cm,height=8cm]{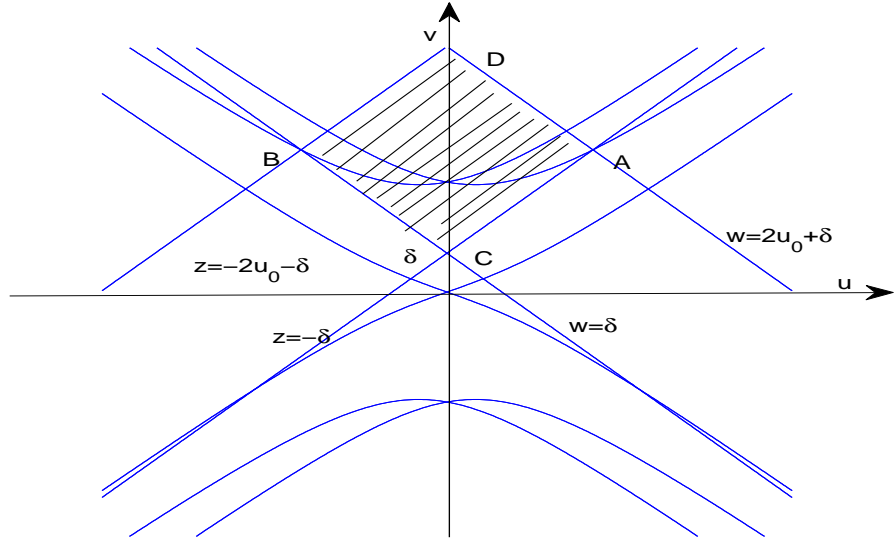}\\
  \caption{Invariant region when $\alpha<-1$}\label{figure5}
  \end{center}
\end{figure}

\begin{Example}\label{Ex1}
For the following special catenoid type surfaces with
 $$E(y)=\left(c\cosh\left(\frac{y}{c}\right)\right)^{\frac{2}{\beta^2-1}}, \quad G(y)=\frac{1}{c^2(\beta^2-1)^2}E(y),$$
  and $$K(y)=-c^2(\beta^2-1)E(y)^{-\beta^2},$$
where $c\neq0, \,\beta\geq \sqrt{2}$ are constants,
one has $E(y)^\prime>0$ whenever $y>0$, and $E(y)^\prime<0$ whenever $y<0$.
All the conditions \eqref{cm1} and \eqref{cm2} for the above invariant region are satisfied when $y<0$.
If we take $\Omega=\{(x,y):x\in \mathbb{R},-y_0\leq y\leq 0\}$,
where $y_0>0$ is arbitrary, then in $\Omega$ the equations (\ref{6})
are parabolic  for $y-$time like. Therefore, when the initial
value $(u(x,-y_0), v(x,-y_0))$  is in the square $ACBD$, the parabolic maximum/minimum
principle ensures that the square  $ACBD$ is an invariant region, 
which yields the $L^\infty$ estimate.
We notice that the surface is just the classical catenoid when $\beta=\sqrt{2}$.
Indeed, 
for the special catenoid type metric, the surface is given by the following function:
$$\mathbf{r}(x,y)=(r_1(x,y),r_2(x,y),r_3(x,y)),$$
with
\begin{eqnarray*}
&&r_1(x,y)=\left(c\cosh\left(\frac{y}{c}\right)\right)^{\frac{1}{\beta^2-1}}\sin(x),\\
&&r_2(x,y)=\left(c\cosh\left(\frac{y}{c}\right)\right)^{\frac{1}{\beta^2-1}}\cos(x),\\
&&r_3(x,y)=\int^y\frac{1}{(\beta^2-1)^2}\left(c\cosh\left(\frac{t}{c}\right)\right)^{\frac{4-2\beta^2}{\beta^2-1}}dt.
\end{eqnarray*}
\end{Example}

\subsection{Helicoid type surfaces:  $E(y)=B(y)^2, F=0, G(y)=1$} \label{helicoid}
For the surfaces with the metric of the form:
\begin{equation}\label{helicoidmetric}
E(y)=B(y)^2,\quad  F=0, \quad G(y)=1, \quad B(y)>0,
\end{equation}
 we can also calculate that
\begin{equation*}
\Gamma^2_{22}=0, \quad \Gamma^2_{12}=0, \quad \Gamma^2_{11}=\frac{-E^\prime}{2E},
\quad \Gamma^1_{22}=0,\quad \Gamma^1_{12}=\frac{E^\prime}{2E},\quad \Gamma^1_{11}=0.
\end{equation*}
and
\begin{equation*}
\begin{split}
\tilde{\Gamma}^2_{22}=\frac{\gamma^\prime}{\gamma}, \qquad \tilde{\Gamma}^2_{12}=0, \qquad \tilde{\Gamma}^2_{11}=\frac{-E^\prime}{2},\\
\tilde{\Gamma}^1_{22}=0,\qquad  \tilde{\Gamma}^1_{12}=
\frac{E^\prime}{2E}+\frac{\gamma^\prime}{2\gamma}, \qquad \tilde{\Gamma}^1_{11}=0.
\end{split}
\end{equation*}
Then
\begin{equation*}
\begin{split}
&f(u,v)=-\frac{E^\prime}{E}u-\frac{E^\prime}{2}(u^2-v^2)u, \\
&g(u,v)=\frac{\gamma^\prime}{\gamma}v-\frac{E^\prime}{2}(u^2-v^2)v.
\end{split}
\end{equation*}

Assuming $$\frac{\gamma^\prime}{\gamma}=a\frac{B'}{B},$$
and using $E(y)=B(y)^2$ and \eqref{K},
 we have $$a\frac{E^\prime}{E}=\frac{K^\prime}{K},$$ that is,
 \begin{equation}\label{hm1}
 K(y)=-k_0E(y)^a,
 \end{equation}
 with some constant $k_0>0$.  Thus,
\begin{equation*}
\begin{split}
&f(u,v)=-\frac{2B^\prime}{B}u-BB^\prime(u^2-v^2)u, \\
&g(u,v)=a\frac{B^\prime}{B}v-BB^\prime(u^2-v^2)v,
\end{split}
\end{equation*}
and
\begin{equation*}
\begin{split}
&f(u,v)+g(u,v)=-\frac{2B^\prime}{B}u+a\frac{B^\prime}{B}v-BB^\prime(u^2-v^2)(u+v),\\
&f(u,v)-g(u,v)=-\frac{2B^\prime}{B}u-a\frac{B^\prime}{B}v-BB^\prime(u^2-v^2)(u-v).
\end{split}
\end{equation*}
Now set
\begin{equation}
\begin{split}
\tilde{u}=Bu, \quad &\tilde{v}=Bv,\\
 \tilde{w}=Bu+Bv=\tilde{u}+\tilde{v},\quad &\tilde{z}=Bu-Bv=\tilde{u}-\tilde{v}.
\end{split}
\end{equation}
Since  $B$ depends only on $y$, we have
\begin{equation}
\begin{split}
&\tilde{w}_y+\lambda_1\tilde{w}_x=\varepsilon \tilde{w}_{xx}-\frac{2\varepsilon v_x\tilde{w}_x}{v}+
(f(u,v)+g(u,v))B+\frac{B^\prime}{B}\tilde{w},\\
&\tilde{z}_y+\lambda_2\tilde{z}_x=\varepsilon \tilde{z}_{xx}-\frac{2\varepsilon v_x\tilde{z}_x}{v}+
(f(u,v)-g(u,v))B+\frac{B^\prime}{B}\tilde{z}.
\end{split}
\end{equation}
Note that
\begin{equation*}
f(u,v)+g(u,v)=-\frac{B^\prime}{B^2}\left(2Bu-aBv+((Bu)^2-(Bv)^2)(Bu+Bv)\right),
\end{equation*}
 then
\begin{equation*}
\begin{split}
R_1(u,v):=&(f(u,v)+g(u,v))B+\frac{B^\prime}{B}\tilde{w}\\
=&-\frac{B^\prime}{B}\left(2\tilde{u}-a\tilde{v}+(\tilde{u}^2-\tilde{v}^2)(\tilde{u}+\tilde{v})\right)
 +\frac{B^\prime}{B}(\tilde{u}+\tilde{v})\\
 =&-\frac{B^\prime}{B}\left(\tilde{u}-(a+1)\tilde{v}+(\tilde{u}^2-\tilde{v}^2)(\tilde{u}+\tilde{v})\right),
\end{split}
\end{equation*}
and similarly,
\begin{equation*}
\begin{split}
R_2(u,v):=&(f(u,v)-g(u,v))B+\frac{B^\prime}{B}\tilde{z}\\
 =&-\frac{B^\prime}{B}\left(\tilde{u}+(a+1)\tilde{v}+(\tilde{u}^2-\tilde{v}^2)(\tilde{u}-\tilde{v})\right).
\end{split}
\end{equation*}
Since $B(y)>0$ is given, it remains to   find the invariant region of $\tilde{w}, \tilde{z}$.
As in the catenoid case in Subsection \ref{catenoid}, we need to analyze the signs of $R_1$ and $R_2$.
 If we assume
 $$B^\prime<0,$$
 or equivalently, from $E=B^2$,
 \begin{equation}\label{hm2}
 E^\prime<0,
 \end{equation}
  and set
\begin{equation*}
\begin{split}
&\psi_1(\tilde{u},\tilde{v}):=\tilde{u}-(a+1)\tilde{v}+(\tilde{u}^2-\tilde{v}^2)(\tilde{u}+\tilde{v}),\\
&\psi_2(\tilde{u},\tilde{v}):=\tilde{u}+(a+1)\tilde{v}+(\tilde{u}^2-\tilde{v}^2)(\tilde{u}-\tilde{v}),
\end{split}
\end{equation*}
then the signs of $R_1(\tilde{u},\tilde{v}),
R_2(\tilde{u},\tilde{v})$ depend on the signs of
$\psi_1(\tilde{u},\tilde{v}), \psi_2(\tilde{u},\tilde{v})$
respectively.
Now we end up with a situation similar to the catenoid case in Subsection \ref{catenoid}
with  $c=1, \, a=2\alpha$, and  the corresponding $\varphi_i$ ($i=1, 2$) in Subsection \ref{catenoid} are
$$\varphi_1(\tilde{u},\tilde{v})=\psi_1(\tilde{u},\tilde{v}), \quad
\varphi_2(\tilde{u},\tilde{v})=\psi_2(\tilde{u},\tilde{v}).$$
We can find the invariant regions just as in catenoid case in Subsection \ref{catenoid}.
Indeed, the invariant region of $(\tilde{u},\tilde{v})$ looks like the same as the invariant region $ACBD$ in Figure \ref{figure3}
when $a=-2$ if we replace the $u-v$ plane by the
$\tilde{u}-\tilde{v}$ plane and  $u_0$ by $\tilde{u}_0=\tilde{u}_0(-2)$ defined below;
and looks like the same as the invariant region $ACBD$
 in Figure \ref{figure5} when $a<-2$ in the $\tilde{u}-\tilde{v}$ plane,
and thus we omit the sketch of the invariant regions.
From the invariant region,
when $a=-2$,
\begin{equation}
-\frac{1-\delta^2}{2\delta}\leq\tilde{u}\leq\frac{1-\delta^2}{2\delta},
\qquad \delta\leq \tilde{v}\leq\frac{1}{\delta};
\end{equation}
and when $a<-2$,
\begin{equation}
-\tilde{u}_0\leq \tilde{u}\leq \tilde{u}_0, \qquad
\delta\leq \tilde{v}\leq 2\tilde{u}_0+\delta,
\end{equation}
where
$$0<\delta<\sqrt{-a-1},$$
 \begin{equation}\label{u03}
 \tilde{u}_0=\tilde{u}_0(a)=\frac{-a-4\delta^2+
\sqrt{a^2-8a\delta^2-16\delta^2}}{8\delta}.
\end{equation}
Note that $0<B(y)\in C^{1,1}(\Omega)$, we can easily
obtain the $L^\infty$ boundedness of $u,v$.

\begin{Example}\label{Ex2}
For the helicoid surface with
$$E(y)=c^2+y^2, \quad F(y)\equiv0, \quad G(y)\equiv1,\quad K(y)=-\frac{c^2}{(c^2+y^2)^2},$$
 where $c\neq 0$, we see that $$B(y)=\sqrt{c^2+y^2},$$ and $B(y)^\prime>0$ for $y>0$,
 and $B(y)^\prime<0$ for $y<0$. If we take $$\Omega=\{(x,y):x\in \mathbb{R},-y_0\leq y\leq 0\},$$
 where $y_0>0$ is an arbitrary constant, then in $\Omega$, the equations in (\ref{6})
are parabolic  for $y-$time like. Therefore, when the  
 initial value $(u(x,-y_0), v(x,-y_0))$ is in the square $ACBD$,  we have the
 invariant region for the solutions. We note that the function for the helicoid in $\mathbb{R}^3$ is
 $\mathbf{r}(x,y)=(y\sin x, y\cos x, cx).$
 \end{Example}

 \begin{Remark}\label{R31}
 We can see from Figures \ref{figure3} and \ref{figure5}  that the $L^\infty$
 estimates also hold for $v<0$ since we can also find the invariant
 regions in the lower half-plane of $u-v$, which are symmetric with the invariant regions in the upper half-plane around the $u$ axis.
 \end{Remark}

\bigskip

\section{$H^{-1}$ Compactness}

In this section  we  shall  prove  the $H^{-1}_{loc}$ compactness
of the approximate viscous solutions.

For the strictly convex entropy $$\eta=\frac{\tilde{M}^2+1}{\tilde{L}},$$ and entropy flux
$$q=\frac{-\tilde{M}^3+\tilde{M}}{\tilde{L}^2},$$
from  the parabolic equations \eqref{gcv}, one has
\begin{equation}\label{10}
\begin{split}
\eta_y+q_x&=\varepsilon \eta_{xx}+\Pi(\tilde{L},\tilde{M})
-2\varepsilon\left(\frac{\tilde{M}^2+1}{\tilde{L}^3}\tilde{L}_x^2-
2\frac{\tilde{M}}{\tilde{L}^2}\tilde{L}_x\tilde{M}_x+
    \frac{\tilde{M}_x^2}{\tilde{L}}\right)\\
&=\varepsilon \eta_{xx}+ \Pi(\tilde{L}, \tilde{M})
-2\varepsilon\left(\frac{\tilde{L}_x^2}{\tilde{L}^3}
+\frac1{\tilde{L}}\left(\frac{\tilde{M}}{\tilde{L}}\tilde{L}_x-\tilde{M}_x\right)^2\right),
\end{split}
\end{equation}
where
\begin{equation*}
\begin{split}
\Pi(\tilde{L},\tilde{M})=&\frac{\tilde{M}^2+1}{\tilde{L}^2}\left(\tilde{\Gamma}^2_{22}\tilde{L}
-2\tilde{\Gamma}^2_{12}\tilde{M}+\tilde{\Gamma}^2_{11}\tilde{N}\right)\\
 &+ \frac{2\tilde{M}}{\tilde{L}}\left(\tilde{\Gamma}^1_{22}\tilde{L}-2\tilde{\Gamma}^1_{12}\tilde{M}
+\tilde{\Gamma}^1_{11}\tilde{N}\right).
\end{split}
\end{equation*}
From $$\tilde{L}=\frac1v, \quad \tilde{M}=-\frac{u}v,$$
we have $$\tilde{L}_x=-\frac{v_x}{v^2}, \quad \tilde{M}_x=\frac{uv_x}{v^2}-\frac{u_x}v,$$
and thus
\begin{equation}\label{10b}
\eta_y+q_x= \varepsilon \eta_{xx}+ \Pi(\tilde{L}, \tilde{M})
-2\varepsilon\left(\frac{v_x^2}{v}+\frac{u_x^2}{v}\right).
\end{equation}
 Due to the $L^\infty$ uniform estimates for $u$ and $v$ in Subsections \ref{catenoid} and \ref{helicoid},
 we have
 \begin{equation}
 0<b_1(\delta) \leq\tilde{L}\leq b_2(\delta), \quad |\tilde{M}|\leq b_3(\delta),
 \end{equation}
 uniformly in $\varepsilon$ in $\Omega$, where $b_1(\delta)$, $b_2(\delta)$
 and $b_3(\delta)$ are positive constants depending on  $\delta>0$. Therefore
 $\Pi(\tilde{L}, \tilde{M})$ is also uniformly bounded in
 $\varepsilon$. Let $$\Omega=\{(x,y):x\in \mathbb{R}, -y_0\leq y\leq0\},$$
 where $y_0>0$ is arbitrary. Choose the test function $\phi\in C^\infty_c(\Omega)$
 satisfying $\phi|_K\equiv1$, $0\leq\phi\leq1$, where $K$ is a compact
 set and $K\subset S=\text{supp }\phi$.  From $\eqref{10b}$,   we have
\begin{equation}
\begin{split}
&\int_{-\infty}^{\infty}\int_{-y_0}^0
2\varepsilon\left(\frac{v_x^2}{v}+\frac{u_x^2}{v}\right)\phi dydx\\
&\leq \int_{-\infty}^{\infty}\int_{-y_0}^0\left(\varepsilon\eta_{xx}-
\eta_{y}-q_x+\Pi(\tilde{L}, \tilde{M})\right)\phi dydx\\
&=\int_{-\infty}^{\infty}\int_{-y_0}^0\left(\varepsilon\eta\phi_{xx}+
\eta\phi_y+q\phi_x+\Pi(\tilde{L}, \tilde{M})\phi\right)dydx\\
&\leq M(\phi)
\end{split}
\end{equation}
for some positive constant $M(\phi)$ uniform in $\varepsilon\in(0,1)$.
Since $v$ is uniformly bounded from below with positive lower bound,
$\varepsilon v_x^2, \, \varepsilon u_x^2$ are bounded in $L^1_{loc}(\Omega)$.
Since $u$ is uniformly bounded and $v$ has uniform positive lower bound, one has
$$\varepsilon\tilde{L}_x^2=\varepsilon\frac{v_x^2}{v^4}\le C\varepsilon v_x^2, \quad
 \varepsilon\tilde{M}_x^2= \varepsilon\left(\frac{uv_x}{v^2}-\frac{u_x}v\right)^2
 \le C\varepsilon(v_x^2+u_x^2),$$
 for some positive constant $C$ uniform in $\varepsilon$,
 then we see that  $\varepsilon \tilde{L}_x^2$ and $\varepsilon \tilde{M}_x^2$ are
uniformly bounded in $L^1_{loc}(\Omega)$.
Noting that $$\varepsilon\tilde{L}_{xx}=\sqrt{\varepsilon}(\sqrt{\varepsilon}\tilde{L}_x)_x,$$
and  for arbitrary $\phi\in C_c^{\infty}(\Omega)$,
\begin{equation}\label{weak}
\begin{split}
&\int_{-\infty}^{\infty}\int_{-y_0}^0\varepsilon \tilde{L}_{xx}\phi dydx=
\int_{-\infty}^{\infty}\int_{-y_0}^0\varepsilon \tilde{L}_x\phi_x dydx\\
&\leq \sqrt{\varepsilon}\left(\int_{\text{supp }\phi}\varepsilon\tilde{L}_x^2 dydx\right)^{\frac{1}{2}}
\left(\int_{\Omega}(\phi_x)^2dydx\right)^{\frac{1}{2}}\\
&\leq C \sqrt{\varepsilon}\rightarrow 0  \text{ as } \varepsilon\rightarrow 0,
\end{split}
\end{equation}
we see that  $\varepsilon\tilde{L}_{xx}$ is compact in $H^{-1}_{loc}(\Omega)$.
From (\ref{gcv}), we have
\begin{equation}
\tilde{M}_x-\tilde{L}_y=\tilde{\Gamma}^2_{22}\tilde{L}-
2\tilde{\Gamma}^2_{12}\tilde{M}+\tilde{\Gamma}^2_{11}\tilde{N}
-\varepsilon\tilde{L}_{xx}.
\end{equation}
 Since  $\tilde{\Gamma}^2_{22}\tilde{L}-2\tilde{\Gamma}^2_{12}\tilde{M}+
 \tilde{\Gamma}^2_{11}\tilde{N}$ is uniformly bounded in $\Omega$,
 it is uniformly bounded in $L^1_{loc}(\Omega)$ and compact in
 $W^{-1,p}_{loc}(\Omega)$ with some $1<p<2$ by the imbedding theorem
 and the Schauder theorem. Therefore $\tilde{M}_x-\tilde{L}_y$ is
 compact in $W^{-1,p}_{loc}(\Omega)$. Moreover,  we  see
 that $\tilde{M}_x-\tilde{L}_y$ is uniformly bounded in
 $W^{-1,\infty}_{loc}(\Omega)$ since  $\tilde{M}$ and $\tilde{L}$ are
 uniformly bounded. Finally, by  Lemma  \ref{lemma1} below,  we conclude  that $\tilde{M}_x-\tilde{L}_y$  is compact in $H^{-1}_{loc}(\Omega)$.
 Similarly, $\tilde{N}_x-\tilde{M}_y$ is also compact in $H^{-1}_{loc}(\Omega)$.
 Since  $\gamma$ is $C^1$, we see that  $M_x-L_y$ and $N_x-M_y$ are  also
  compact in $H^{-1}_{loc}(\Omega).$

We record the following useful lemma (see \cite{Chen1, Ta1}) here:
\begin{Lemma}\label{lemma1}
Let $\Omega\in\mathbb{R}^n$ be a open set, then
(compact set of $W^{-1,q}_{loc}(\Omega)$ ) $\bigcap$
(bounded set of $W_{loc}^{-1,r}(\Omega)$ )
$\subset$ (compact set of $H^{-1}_{loc}(\Omega)$).
where $q$ and $r$ are constants, $1<q\leq2<r.$
\end{Lemma}

\bigskip

\section{Main Theorems and Proofs}\label{S5}

In this section, we shall state our main results and also give the proof.

 In the previous sections, we have established the  $L^{\infty}$ uniform estimate and $H_{loc}^{-1}$
 compactness  of the viscous approximate solutions to \eqref{gcv} for some special metrics of the form
  $$ds^2=Edx^2+Gdy^2.$$
 The corresponding Gauss curvature $K(x,y)=K$ has
  the following form (see \cite{HH}):
\begin{equation}\label{curvature}
K=\frac{1}{4EG}\left(\frac{E_y^2+E_xG_x}{E}+\frac{G_x^2+E_yG_y}{G_y}-2(E_{yy}+G_{xx})\right).
\end{equation}
To prove the existence of isometric immersion, first let us recall the
following compensated compactness framework in Theorem 4.1 of \cite{CSW}:

\begin{Lemma}\label{lemma2}
Let a sequence of functions $(L^\varepsilon, M^\varepsilon, N^\varepsilon)(x,y)$,
defined on an open subset $\Omega\subset\mathbb{R}^2$, satisfy the
following framework:
\begin{itemize}
\item[(W.1)]$(L^\varepsilon, M^\varepsilon, N^\varepsilon)(x,y)$ is uniformly
bounded almost everywhere in $\Omega\subset\mathbb{R}^2$ with respect to $\varepsilon$;
\item[(W.2)]$M^\varepsilon_x-L^\varepsilon_y$ and $N^\varepsilon_x-M^\varepsilon_y$ are
compact in $H^{-1}_{loc}(\Omega)$;
\item[(W.3)]There exist $o_j^\varepsilon(1)$, $j=1,2,3$, with $o_j^\varepsilon(1)\rightarrow0$
in the sense of distributions as $\varepsilon\rightarrow0$ such that
\begin{equation}\label{c}
\begin{split}
 &M_x^\varepsilon-L_y^\varepsilon=\Gamma^2_{22}L^\varepsilon-2\Gamma^2_{12}M^\varepsilon
 +\Gamma^2_{11}N^\varepsilon+o_1^\varepsilon(1),\\
 &N_x^\varepsilon-M_y^\varepsilon=-\Gamma^1_{22}L^\varepsilon+2\Gamma^1_{12}M^\varepsilon
 -\Gamma^1_{11}N^\varepsilon+o_2^\varepsilon(1),\\
 \end{split}
\end{equation}
and
\begin{equation}\label{G}
L^\varepsilon N^\varepsilon-(M^\varepsilon)^2=K+o_3^\varepsilon(1).
\end{equation}
\end{itemize}
Then there exists a subsequence (still labeled) $(L^\varepsilon, M^\varepsilon, N^\varepsilon)$
converging weak-star in $L^\infty$ to $(L, M, N)(x,y)$ as $\varepsilon\rightarrow0$ such that
\begin{itemize}
\item[(1)]$ (L, M, N)$ is also bounded in $\Omega\subset\mathbb{R}^2$;
\item[(2)]the Gauss equation (\ref{G}) is weakly continuous with respect to the
 subsequence $(L^\varepsilon, M^\varepsilon, N^\varepsilon)$ converging weak-star in
 $L^\infty$ to $(L, M, N)(x,y)$ as $\varepsilon\rightarrow0$;
\item[(3)] The Codazzi equations (\ref{c})  as $\varepsilon\to 0$ hold for $(L, M, N)$ in the sense
of distribution.
\end{itemize}
\end{Lemma}

We now present the results on the existence of isometric immersion of surfaces with the two types of metrics studied in Section 3.

\begin{Definition}\label{def1}
A Riemannian metric on a two-dimensional manifold is called a {\em catenoid  metric}  if it is of the form
$$ds^2=E(y)dx^2+cE(y)dy^2,$$
with  $$c>0, \quad E(y)>0, \quad E(y)^\prime<0 \text{ for } y<0,$$  and
the corresponding  Gauss curvature is of the form $$K(y)=-k_0E(y)^{-\beta^2}$$ with constants $k_0>0$ and $\beta\geq\sqrt{2}$.
\end{Definition}

 From  the Definition \ref{def1} and  the formula in (\ref{curvature}), we see that $E(y)$ satisfies
the following ordinary differential equation:
\begin{equation}\label{13}
(E(y)^\prime)^2-EE(y)^{\prime\prime}=-2k_0E(y)^{2-\beta^2}.
\end{equation}
We can solve it through the following process utilizing the method of \cite{PZ}.
Set $$E(y)=e^{w(y)},$$  then
$$E(y)^\prime=e^{w(y)}w(y)^\prime, \quad E(y)^{\prime\prime}=e^{w(y)}\left((w(y)^\prime)^2+w(y)^{\prime\prime}\right),$$
and (\ref{13}) becomes
\begin{equation}
w(y)^{\prime\prime}=2k_0e^{-\beta^2w(y)}.
\end{equation}
Let $f(w(y))=w(y)^\prime.$ After differentiating respect to $y$, we get
$$f(w)^\prime w(y)^\prime=w(y)^{\prime\prime},$$
i.e.,
$$f(w)f(w)^{\prime}=2k_0e^{-\beta^2w},$$
therefore, noting that $E'<0$ implies $w'<0$,
$$f(w)=\frac{d w}{d y}=-\sqrt{C_1-\frac{4k_0}{\beta^2} e^{-\beta^2w}}.$$
Then one has
$$y=C_2-\int\frac{dw}{\sqrt{C_1-\frac{4k_0}{\beta^2}e^{-\beta^2w}}},$$
Denote the right side of the above equation by $h(w)$,
then $w(y)=h^{-1}(y)$ and $E(y)=e^{h^{-1}(y)}$,
where $h^{-1}(y)$ is the inverse function of $h(w)$, $C_1$ and $C_2$ depend on the
the value of $w(0)$ and $w(0)^{\prime}$.
Then the catenoid metric is  of the form:
\begin{equation}
ds^2=e^{h^{-1}(y)}dx^2+ce^{h^{-1}(y)}dy^2.
\end{equation}
We note that the special catenoid type metric given in Example \ref{Ex1}:  
 $$E(y)=\left(c\cosh\left(\frac{y}{c}\right)\right)^{\frac{2}{\beta^2-1}}, \quad
 G(y)=\frac{1}{c^2(\beta^2-1)^2}E(y),$$
with  $K(y)=-c^2(\beta^2-1)E(y)^{-\beta^2}$
is a catenoid metric in the sense of Definition \ref{def1}.

\begin{Definition}\label{def2}
A Riemannian metric on a two-dimensional manifold is called a {\em helicoid metric}
if it is of the form $$ds^2=E(y)dx^2+dy^2,$$
with  $$E(y)>0, \quad E(y)^\prime<0 \text{ for } y<0,$$
and the corresponding Gauss curvature is of the form $$K(y)=-k_0E(y)^a$$ with constants $k_0>0$ and $a\leq-2$.
\end{Definition}

From Definition \ref{def2} and the formula (\ref{curvature}), $E(y)$ satisfies  the following
ordinary differential equation:
\begin{equation}\label{14}
(E(y)^\prime)^2-2EE(y)^{\prime\prime}=-4k_0E(y)^{2+a}.
\end{equation}
Set  $$E(y)=w(y)^2,$$  then (\ref{14}) becomes
$$w(y)^{\prime\prime}=k_0w(y)^{2a+1}.$$
Letting $g(w(y))=w(y)^\prime$, and differentiating respect to $y$, one has
$$w(y)^{\prime\prime}=u(w)^\prime w(y)^{\prime},$$
i.e.,
$$g(w)^{\prime}g(w)=k_0w^{2a+1}.$$
Thus
$$g(w)=\frac{d w}{d y}=-\sqrt{C_1+\frac{k_0}{a+1}w^{2a+2}}.$$
Then
$$y=C_2-\int\frac{dw}{\sqrt{{C_1+\frac{k_0}{a+1}}w^{2a+2}}}.$$
Denote the right side of the above equation by $h(w)$, then $w(y)=h^{-1}(y),$ and
$E(y)=(h^{-1}(y))^2.$ Therefore, the helicoid metric is
$$ds^2=(h^{-1}(y))^2dx^2+dy^2,$$
where $h^{-1}(y)$ is the inverse function of
$$h(w)=C_2-\int\frac{dw}{\sqrt{C_1+\frac{k_0}{a+1}w^{2a+2}}}.$$
Similar to the  catenoid metric,
$C_1$ and $C_2$ depend on the value of $w(0)$ and $w(0)^{\prime}$.
As an example, the helicoid surface with  $$E(y)=c^2+y^2,\quad K(y)=-\frac{c^2}{(c^2+y^2)^2}$$
is a two-dimensional Riemannian manifold with the helicoid metric
in the sense of  Definition \ref{def2}.

Now we prove the existence of $C^{1,1}$ isometric immersion of surfaces with the  above two types of metrics
 into $\mathbb{R}^3$ by  using   Lemma \ref{lemma2}.

\begin{Theorem}\label{Th1}
For any given $y_0>0$, let the initial data
\begin{equation}\label{ini}
(u, v)|_{y=-y_0}=(\bar{u}_0(x), \bar{v}_0(x)):=(u(x,-y_0), v(x, -y_0))  
\end{equation}
satisfy the following conditions:
$$\bar{u}_0+\bar{v}_0 \text{ and } \bar{u}_0-\bar{v}_0 \text{ are bounded},$$
and
$$\inf_{x\in \mathbb{R}}\left(\bar{u}_0+\bar{v}_0\right)>0, \quad \sup_{x\in \mathbb{R}}\left(\bar{u}_0-\bar{v}_0\right)<0.$$
Then,  for the catenoid metric in the sense of Definition \ref{def1},
the Gauss-Codazzi system (\ref{Gc}) has a weak solution in
$\Omega=\{(x,y): x\in \mathbb{R}, -y_0\leq y\leq0\}$ with the initial data (\ref{ini}).
\end{Theorem}

\begin{Remark}\label{R51}
For the initial data $(\bar{u}_0(x), \bar{v}_0(x))$ satisfying the conditions in Theorem \ref{Th1}, there exists a constant
$\delta>0$, such that
\begin{equation*}
\begin{split}
(\bar{u}_0(x), \bar{v}_0(x))\in ACBD:=\{(u,v):  \;\;  &\delta\le w=u+v\le 2u_0+\delta, \\
                                                                        &-(2u_0+\delta)\le z=u-v\le -\delta\},
\end{split}
\end{equation*}
for the catenoid type metrics, where $u_0$ is defined in \eqref{u02}. That is, $(\bar{u}_0(x), \bar{v}_0(x))$  lies in the invariant region $ACBD$ sketched in Figure 3 or Figure 5, i.e.,  $\bar{u}_0(x)$ is bounded,
and $\bar{v}_0(x)$ has positive lower bound and  upper bound.
\end{Remark}

\begin{proof}
First for the initial data \eqref{ini}, the corresponding initial data for $\tilde{L}, \tilde{M}, \tilde{N}$ is
 $$\tilde{L}_0(x)=\frac{1}{\bar{v}_0(x)}, \quad \tilde{M}_0(x)=\frac{-\bar{u}_0(x)}{\bar{v}_0(x)}, \quad
\tilde{N}_0(x)=\frac{\bar{u}_0(x)^2-\bar{v}_0(x)^2}{\bar{v}_0(x)}.$$
We  use system (\ref{6}) to obtain
the approximate viscous solutions and their $L^\infty$ estimate, and then use (\ref{gcv}) to obtain the $H^{-1}_{loc}$ compactness.

\emph{Step 1}. We mollify  the initial data \eqref{ini} as
\begin{equation*}
\bar{u}_0^\varepsilon(x)=\bar{u}_0(x)\ast j^\varepsilon, \quad \bar{v}_0^\varepsilon(x)=\bar{v}_0(x)\ast j^\varepsilon,
\end{equation*}
where $j^\varepsilon$ is  the standard mollifier and $\ast$ is for the convolution.
From the above Remark \ref{R51}, there exists a $\delta>0$, such that,
\begin{equation*}
|\bar{u}_0^\varepsilon(x)|\leq M(\delta) , \quad 0<\delta\leq \bar{v}_0^\varepsilon(x)\leq M(\delta),
\end{equation*}
where $M(\delta)$ is positive constant depending on $\delta$, and
\begin{equation*}
(\bar{u}_0^\varepsilon(x), \bar{v}_0^\varepsilon(x))\rightarrow (\bar{u}_0(x), \bar{v}_0(x)) \text{ as } \varepsilon\rightarrow 0,  \qquad a.e.
\end{equation*}

\emph{Step 2}. The local existence of (\ref{gcv}) can be obtained by the standard theory, hence we have the
local existence of (\ref{6}).  From Section 3,
we have the $L^\infty$ estimate for the approximate solution of (\ref{6}),
 then we can obtain the global existence of the approximate solution as follows.
We observe that  the first equation of (\ref{6}) is of
the divergence form, thus  the estimate of $u_x$ can be handled in the standard way.
However, the second equation of \eqref{6} is not in divergence form.
By differentiation of the second equation with respect to $x$, we get
\begin{equation*}
(v_x)_y+(uv_x-vu_x)_x=g(u,v)_x+\varepsilon (v_x)_{xx}-\Big(\frac{2\varepsilon v_x^2}{v}\Big)_x,
\end{equation*}
where we omit $\varepsilon$ and still use $u(x,y), v(x,y)$   for the approximate solution of (\ref{6}).
Let $G(x,y)$ is the heat kernel of $h_y=\varepsilon h_{xx}$, then we can solve the above
equation as the following:
\begin{equation*}
\begin{split}
v_x=&\int_{-\infty}^{\infty}G(x-\xi,y)\bar{v}_{0x}^\varepsilon(\xi)d\xi\\
 &+\int_{-y_0}^{0}\int_{-\infty}^{\infty}\left(g(u,v)_x-(uv_x-vu_x)_x
 +\Big(\frac{2\varepsilon v_x^2}{v}\Big)_x\right)G(x-\xi,y-\eta) d\xi d\eta \\
 =&\int_{-\infty}^{\infty}G(x-\xi,y)\bar{v}_{0x}^\varepsilon(\xi)d\xi\\
 &+\int_{-y_0}^{0}\int_{-\infty}^{\infty}\left(-g(u,v)+(uv_x-vu_x)
 +\frac{2\varepsilon v_x^2}{v}\right)G(x-\xi,y-\eta)_x d\xi d\eta.
\end{split}
\end{equation*}
Therefore
\begin{equation*}
\begin{split}
&||v_x||_{C^0}\leq  \int_{-\infty}^{\infty}G(x-\xi,y)|\bar{v}_{0x}^\varepsilon(\xi)|d\xi\\
 &\qquad\qquad +\int_{-y_0}^{0}\int_{-\infty}^{\infty}\left(|g(u,v)|+|(uv_x-vu_x)|+
 \left|\frac{2\varepsilon v_x^2}{v}\right|\right)G(x-\xi,y-\eta)_x| d\xi d\eta\\
&\leq  ||\bar{v}_{0x}^\varepsilon||_{C^0}+M(\delta)(||u_x||_{C^0}|+||v_x||_{C^0}+\varepsilon||v_x||_{C^0})
\int_{-y_0}^{0}\int_{-\infty}^{\infty}|G(x-\xi,y-\eta)_x| d\xi d\eta,
\end{split}
\end{equation*}
where $||\cdot||_{C^0}$ stands for the norm of continuous functions. The above
process yields the $C^0$ estimate of $v_x$.
Similarly,  we can estimate $C^k$ (for any integer $k\ge 1$) norms
of $u(x,y), v(x,y)$,  which  can be bounded by the $C^k$ norms
 of $\bar{u}_0^\varepsilon$ and $\bar{v}_0^\varepsilon$.  Thus, we obtain   the global existence
  of smooth solutions to the system (\ref{6}).

\emph{Step 3.}  In Section 3 we have proved $L^\infty$ boundness of
$\tilde{L}, \tilde{M}, \tilde{N}$,  and then $L, M, N$ is in $L^\infty$
for $\gamma\in C^1$. By the reverse process of Section 2, we can reformulate
the equations (\ref{6}) of $u$ and $v$ as the equations (\ref{gcv}).
Therefore, as in Section 4 we also obtain the $H^{-1}_{loc}$ compactness.
 So we have proved that our approximate solutions satisfy (W.1)
  and (W.2)  in the framework of Lemma \ref{lemma2}.
Furthermore,  from Section 2, we have
\begin{equation*}
\begin{split}
&(\tilde{M}_x-\tilde{L}_y)-(\tilde{\Gamma}^2_{22}\tilde{L}-
2\tilde{\Gamma}^2_{12}\tilde{M}+\tilde{\Gamma}^2_{11}\tilde{N})
=-\varepsilon\tilde{L}_{xx},\\
&(\tilde{N}_x-\tilde{M}_y)-(-\tilde{\Gamma}^1_{22}\tilde{L}+
2\tilde{\Gamma}^1_{12}\tilde{M}-\tilde{\Gamma}^1_{11}\tilde{N})
=-\varepsilon\tilde{M}_{xx}.
\end{split}
\end{equation*}
As in  Section 4, from (\ref{weak}), we can get that, as $\varepsilon\to 0$,
$$\varepsilon\tilde{L}_{xx}\rightarrow 0,$$
in the sense of distribution,
and then
$$\tilde{M}_x-\tilde{L}_y=\tilde{\Gamma}^2_{22}\tilde{L}-
2\tilde{\Gamma}^2_{12}\tilde{M}+\tilde{\Gamma}^2_{11}\tilde{N}+o_1(1)$$
holds in the sense of distribution. Similarly,
$$\tilde{N}_x-\tilde{M}_y=-\tilde{\Gamma}^1_{22}\tilde{L}+
2\tilde{\Gamma}^1_{12}\tilde{M}-\tilde{\Gamma}^1_{11}\tilde{N}+o_2(1)$$
also holds in the sense of distribution.
Here $o_1(1), o_2(1) \to 0$ as $\varepsilon\to 0$.
We note that the Gauss equation holds exactly for the viscous approximate solutions.
Therefore (W.3) is satisfied.
Consequently, we complete the proof of the theorem and obtain the isometric immersion  of the surface with the
 catenoid metric in $\mathbb{R}^3$ using Lemma \ref{lemma2}.
\end{proof}

Since we also obtained the $L^\infty$ estimate and the $H^{-1}_{loc}$ compactness for the helicoid metric in the
previous sections,  we can have the isometric immersion of surfaces with the helicoid metric  in $\mathbb{R}^3$ just as the case for the catenoid metric.

\begin{Theorem}\label{Th2}
For any given $y_0>0$, let the initial data
\begin{equation*}
(u, v)|_{y=-y_0}=(\bar{u}_0(x), \bar{v}_0(x)):=(u(x,-y_0), v(x, -y_0))  
\end{equation*}
satisfy the following conditions:
$$\bar{u}_0+\bar{v}_0 \text{ and } \bar{u}_0-\bar{v}_0 \text{ are bounded},$$
and
$$\inf_{x\in \mathbb{R}}\left(\bar{u}_0+\bar{v}_0\right)>0, \quad \sup_{x\in \mathbb{R}}\left(\bar{u}_0-\bar{v}_0\right)<0.$$
Then for the helicoid metric in the sense of Definition \ref{def2},
the Gauss-Codazzi system (\ref{Gc}) has a weak solution in
$\Omega=\{(x,y): x\in \mathbb{R}, -y_0\leq y\leq0\}$ with the initial data (\ref{ini}).
%
\end{Theorem}

\begin{Remark}
 Although the catenoid metric for $\beta>\sqrt{2}$ has also
been studied by Chen-Slemrod-Wang in \cite{CSW}, their
$L^\infty$ estimate is different from ours, especially for $L$.  
In addition, we also  prove
the isometric immersion of catenoid for $\beta=\sqrt{2}$.
\end{Remark}

\begin{Remark}
By Remark \ref{R31}, if we replace the second condition on the initial data in Theorem \ref{Th1} and Theorem \ref{Th2} by the following:
$$\sup_{x\in \mathbb{R}}\left(\bar{u}_0+\bar{v}_0\right)<0, \quad \inf_{x\in \mathbb{R}}\left(\bar{u}_0-\bar{v}_0\right)>0,$$
then both theorems still hold.
\end{Remark}

\begin{Remark}
By even symmetry from the weak solution in $\Omega=\{(x,y): x\in \mathbb{R}, \, -y_0\leq y\leq0\}$ obtained in Theorem \ref{Th1} and Theorem \ref{Th2}, we can obtain the weak solution in $\Omega'=\{(x,y): x\in \mathbb{R}, \, 0\le y\le y_0\}$.
\end{Remark}

\bigskip

\section*{Acknowledgments}
F. Huang's research was supported in part  by NSFC Grant No. 11371349,
National Basic Research Program of China (973 Program) under Grant
No. 2011CB808002, and the CAS Program for Cross $\&$ Cooperative
Team of the Science $\&$ Technology Innovation.
D. Wang's research was supported in part by the NSF Grant DMS-1312800 and NSFC Grant No. 11328102.

\end{document}